\documentclass{amsart}
\usepackage{amssymb,amsmath,amsthm}

\theoremstyle{plain}
\newtheorem{lemma}{Lemma}[section]
\newtheorem{theorem}[lemma]{Theorem}
\newtheorem{corollary}[lemma]{Corollary}
\newtheorem{sublemma}[lemma]{Sublemma}

\theoremstyle{definition}
\newtheorem{example}[lemma]{Example}

\numberwithin{equation}{section}

\newcommand{\K}{\mathbf{k}}
\newcommand{\f}{\varphi}
\newcommand{\tp}{\overline}

\DeclareMathOperator{\gr}{gr}
\DeclareMathOperator{\ak}{AK}
\DeclareMathOperator{\trdeg}{trdeg}
\DeclareMathOperator{\grdeg}{grdeg}
\DeclareMathOperator{\Spec}{Spec}
\DeclareMathOperator{\EXP}{EXP}
\DeclareMathOperator{\ch}{char}
\DeclareMathOperator{\Frac}{Frac}
\DeclareMathOperator{\Aut}{Aut}
\begin{document}
\title{On automorphisms of Danielewski surfaces}
\author{Anthony J. Crachiola}
\address{Department of Mathematics\\
Wayne State University\\
Detroit, MI 48202\\
USA} \email{crach@math.wayne.edu}

\date{June 18, 2004}

\begin{abstract}
We develop techniques for computing the AK invariant of a domain
with arbitrary characteristic. We use these techniques to describe
for any field $\K$ the automorphism group of $\K[X,Y,Z]/(X^n Y -
Z^2 - h(X)Z)$, where $h(0) \ne 0$ and $n \geq 2$, as well as the
isomorphism classes of these algebras.
\end{abstract}
\maketitle

\section{Introduction}

All rings in this paper are commutative with identity. Throughout
this paper, let $\K$ denote  field of arbitrary characteristic,
and let $\K^* = \K \setminus 0$. For a ring $A$, let $A^{[n]}$
denote the polynomial ring in $n$ indeterminates over $A$. Let
$\mathbf{C}$ denote the field of complex numbers.

One fundamental algebraic problem is to describe the automorphism
group of a given algebra. From the perspective of algebraic
geometry, this means describing automorphisms of a given affine
variety. The Jung-van der Kulk \cite{jung,kulk} theorem provides
the answer for the polynomial ring $\K^{[2]}$, i.e. the affine
plane $\K^2$, but for higher dimensional affine spaces the problem
is still open. For other affine varieties there is no general
approach to solving the problem.

Let $A$ be an algebra with characteristic zero. A new tool
appeared in the 1990s when Leonid Makar-Limanov \cite{ML1}
introduced the AK invariant (more commonly known now as the
Makar-Limanov invariant) of $A$ as the intersection of the kernels
of all locally nilpotent derivations on $A$. Each automorphism of
$A$ restricts to an automorphism of the subalgebra $\ak(A)$,
making this invariant useful in describing the automorphism group
of $A$. As one demonstration Makar-Limanov has computed the
automorphism group of a surface $x^n y = P(z)$ over $\mathbf{C}$
\cite{ML2}. The successful application of the AK invariant to this
and other algebro-geometric problems, such as the linearization
conjecture for $\mathbf{C}^*$-actions on $\mathbf{C}^3$
\cite{KKMLR}, has contributed to its current popular status in
algebraic geometry.

The AK invariant as defined by Makar-Limanov loses its potency for
algebras with prime characteristic $p$ because the kernel of each
derivation becomes much larger, containing the $p$th power of
every element. To the author's knowledge, while the AK invariant
is still bearing fruit, all the research is being conducted under
the restriction of zero characteristic. Now, in the characteristic
zero arena locally nilpotent derivations on an algebra $A$ are
interchangeable with algebraic additive group actions on
$\Spec(A)$, and unlike derivations these actions maintain their
attractive properties for prime characteristic algebras. We can
use this point of view to redefine the AK invariant for rings with
arbitrary characteristic.

In the first part of this paper we explain how to use the AK
invariant for domains of arbitrary characteristic. While complex
algebraic geometry can utilize the topological properties of the
complex numbers, the techniques in this paper rely only on
algebraic structures and do not even require algebraic closure of
the ground field. In fact, the results in this paper are valid
over any field of any characteristic. We next compute the AK
invariant of the algebra $R = \K[X,Y,Z]/(X^n Y - Z^2 - h(X)Z)$,
where $n \geq 2$ and $h(0) \ne 0$, and use it to describe the
automorphism group of $R$.

Over $\mathbf{C}$ the algebra $R$ is the coordinate ring of the
surface $x^n y = z^2 + h(x)z$. This is a generalization of the
celebrated Danielewski surface which plays a role in the
cancellation problem for affine varieties. Here is the one form of
the problem. If $V$ and $W$ are affine varieties over $\K$, does
$V \times \K^n \cong W \times \K^n$ imply $V \cong W$? If $\dim(V)
= \dim(W) =1$ the answer is affirmative. This was shown
algebraically by Shreeram Abhyankar, Paul Eakin, and William
Heinzer \cite{AEH}. (As a side remark, there is a new proof of
this fact which employs the AK invariant \cite{CML}.) For an
algebraist the problem is reformulated as follows. If $A$ and $B$
are $\K$-algebras, does $A^{[n]} \cong B^{[n]}$ imply $A \cong B$?
Mel Hochster published the first counterexample \cite{hochster} in
1972, the same year of the Abhyankar-Eakin-Heinzer paper, using
4-dimensional algebras over the field of real numbers. For a
counterexample over an algebraically closed field the world waited
until 1989. It is due to Wlodzimierz Danielewski
\cite{danielewski} who never published the result. For a published
treatment on the Danielewski surfaces, refer to the paper of
Karl-Heinz Fieseler \cite{fieseler} which gives a classification
of normal surfaces equipped with a nontrivial
$\mathbf{C}^+$-action. Here is Danielewski's original
counterexample. Let $V$ and $W$ be surfaces over $\mathbf{C}$
given by $xy = z^2 + z$ and $x^2 y = z^2 + z$, respectively.
Danielewski showed geometrically that the cylinders $V \times
\mathbf{C}$ and $W \times \mathbf{C}$ are isomorphic while $V$ and
$W$ are not. To explain the isomorphism of cylinders, Danielewski
showed that $V$ and $W$ are total spaces of some principal
$\mathbf{C}^+$-bundles over a line with a double point
$\widetilde{\mathbf{C}}$, and that each of these total spaces is
isomorphic to $V \times_{\widetilde{\mathbf{C}}} W$ which is a
trivial bundle over each $V$ and $W$. To distinguish $V$ and $W$
Danielewski used the first homology group at infinity. The AK
invariant can also tell them apart.

Now let $V_i$ be the surface in $\mathbf{C}^3$ given by $x^{n_i} y
= z^2 + h_i(x)z$, $i=1,2$, where $h_i(0) \ne 0$. In the same
unpublished paper \cite{danielewski} Danielewski conjectured that
$V_1 \cong V_2$ if and only if $n_1 = n_2$ and $h_2(x) = \lambda
h_1(\mu x)$ for some $\lambda, \mu \in \mathbf{C}^*$. J\"{o}rn
Wilkens \cite{wilkens} proved this conjecture and also that, as
with Danielewski's original example, the cylinders over $V_1$ and
$V_2$ are isomorphic for any $n_1,n_2$ and any $h_1(x),h_2(x)$. At
the end of our paper we shall revisit this connection to
cancellation as an application of the AK invariant. We describe
the isomorphism classes of $x^n y = z^2 + h(x)z$ over any field
$\K$ and show how an isomorphism of cylinders can be explained
algebraically.

Danielewski surfaces continue to be a source of interest for
current research. In addition to this paper, see for instance
\cite{dubouloz,freudenburg,shpilrain}.

The main inspiration for this paper is the paper \cite{ML2} of
Leonid Makar-Limanov in which similar results are achieved on the
surface $x^n y = P(z)$ over $\mathbf{C}$.

\section{Methods}

\subsection*{Exponential maps, the AK invariant, and locally finite iterative higher derivations}

Let $A$ be a $\K$-algebra. Suppose $\f:A \to A^{[1]}$ is a
$\K$-algebra homomorphism. We write $\f = \f_U:A \to A[U]$ if we
wish to emphasize an indeterminate $U$. We say that $\f$ is an
\emph{exponential map on $A$\/} if it satisfies the following two
additional properties.
\begin{itemize}
\item[(i)] $\varepsilon_0 \f_U$ is the identity on $A$, where
$\varepsilon_0:A[U] \to A$ is evaluation at $U=0$.\label{exp1}

\item[(ii)] $\f_S \f_U = \f_{S+U}$, where $\f_S$ is extended by
$\f_S(U)=U$ to a homomorphism $A[U] \to A[S,U]$.\label{exp2}
\end{itemize}
(When $A$ is the coordinate ring of an affine variety $\Spec(A)$
over $\K$, the exponential maps on $A$ correspond to algebraic
actions of the additive group $\K^+$ on $\Spec(A)$~\cite[\S
9.5]{essen}.)

Given an exponential map $\f: A \to A[U]$, set $\f(U) = U$ to
obtain an automorphism of $A[U]$ with inverse $\f_{-U}$. Consider
the map $\varepsilon_1 \f :A \to A$, where $\varepsilon_1: A[U]
\to A$ is evaluation at $U=1$. One can check that $\varepsilon_1
\f$ is an automorphism of $A$ with inverse $\varepsilon_1
\f_{-U}$.

Define
\[
A^{\f} = \{ a \in A \,|\, \f(a)=a \},
\]
a subalgebra of $A$ called the \emph{ring of $\f$-invariants}. Let
$\EXP(A)$ denote the set of all exponential maps on $A$. We define
the \emph{AK invariant}, or \emph{ring of absolute constants of
$A$\/} as
\[
\ak(A) = \bigcap_{\f \in \EXP(A)} A^{\f}.
\]
This is a subalgebra of $A$ which is preserved by isomorphism.
Indeed, any isomorphism $f:A \to B$ of $\K$-algebras restricts to
an isomorphism $f: \ak(A) \to \ak(B)$. To understand this, observe
that if $\f \in \EXP(A)$ then $f \f f^{-1} \in \EXP(B)$. Remark
that $\ak(A) = A$ if and only if the only exponential map on $A$
is the standard inclusion $\f(a) = a$ for all $a \in A$.

It is often helpful to view a given $\f \in \EXP(A)$ as a sequence
in the following way. For each $a \in A$ and each natural number
$n$, let $D^n(a)$ denote the $U^n$-coefficient of $\f(a)$. Let $D
= \{ D^0, D^1, D^2, \ldots \}$. To say that $\f$ is a $\K$-algebra
homomorphism is equivalent to saying that the sequence
$\{D^i(a)\}$ has finitely many nonzero elements for each $a \in
A$, that $D^n:A \to A$ is $\K$-linear for each natural number $n$,
and that the Leibniz rule
\[
D^n(ab) = \sum_{i+j=n} D^i(a) D^j(b)
\]
holds for all natural numbers $n$ and all $a,b \in A$. The above
properties (i) and (ii) of the exponential map $\f$ translate into
the following properties of $D$.
\begin{itemize}
\item[(i')] $D^0$ is the identity map.

\item[(ii')] (iterative property) For all natural numbers $i,j$,
\[
D^i D^j = \binom{i+j}{i} D^{i+j}.
\]
\end{itemize}

Due to all of these properties, the collection $D$ is called a
\emph{locally finite iterative higher derivation on $A$}. More
generally, a \emph{higher derivation on $A$\/} is a collection $D
= \{D^i\}$ of $\K$-linear maps on $A$ such that $D^0$ is the
identity and the above Leibniz rule holds. The notion of higher
derivations is due to H.~Hasse and F.K.~Schmidt~\cite{hasse}. When
the characteristic of $A$ is zero, each $D^i$ is determined by
$D^1$, which is a locally nilpotent derivation on $A$. In this
case, $\f = \exp(U D^1) = \sum_i \frac{1}{i!}(U D^1)^i$ and
$A^{\f}$ is the kernel of $D^1$. So we retrieve the original
characteristic zero definition of $\ak(A)$ given by L.
Makar-Limanov as the intersection of the kernels of locally
nilpotent derivations on $A$.

The above discussion of exponential maps, locally finite iterative
higher derivations, and the AK invariant makes sense more
generally for any (not necessarily commutative) ring. However, we
will not require this generality.

\subsection*{Degree functions and related lemmas}

Given an exponential map $\f: A \to A[U]$ on a domain $A$ over
$\K$, we can define the \emph{$\f$-degree\/} of an element $a \in
A$ by $\deg_{\f}(a) = \deg_U(\f(a))$ (where $\deg_U(0) = -
\infty$). Note that $A^{\f}$ consists of all elements of $A$ with
non-positive $\f$-degree. The function $\deg_{\f}$ is a degree
function on $A$ because it satisfies the following two properties
for all $a,b \in A$.
\begin{itemize}
\item[(i)] $\deg_{\f}(ab) = \deg_{\f}(a) + \deg_{\f}(b)$.

\item[(ii)] $\deg_{\f}(a+b) \leq \max \{\deg_{\f}(a),\deg_{\f}(b)
\}$.
\end{itemize}
Equipped with these notions, we now collect some useful facts.

\begin{lemma}\label{L:facts}
Let $\f$ be an exponential map on a domain $A$ over $\K$. Let $D =
\{D^i\}$ be the locally finite iterative higher derivation
associated to $\f$.

\renewcommand{\theenumi}{\alph{enumi}}
\begin{itemize}
\item[(a)] If $a,b \in A$ such that $ab \in A^{\f} \setminus 0$,
then $a,b \in A^{\f}$. In other words, $A^{\f}$ is factorially
closed in $A$.\label{factsa}

\item[(b)] $A^{\f}$ is algebraically closed in $A$.\label{factsb}

\item[(c)] For each $a \in A$, $\deg_{\f}(D^i(a)) \leq
\deg_{\f}(a) - i$. In particular, if $a \in A \setminus 0$ and $n
= \deg_{\f}(a)$, then $D^n(a) \in A^{\f}$.\label{factsc}
\end{itemize}
\end{lemma}

\begin{proof}
(a): We have $0 = \deg_{\f}(ab) = \deg_{\f}(a) + \deg_{\f}(b)$,
which implies that
$\deg_{\f}(a) = \deg_{\f}(b) = 0$.\\
(b): If $a \in A \setminus 0$ and $c_n a^n + \cdots + c_1 a + c_0
= 0$ is a polynomial relation with minimal possible degree $n \geq
1$, where each $c_i \in A^{\f}$ with $c_0 \ne 0$, then $a(c_n
a^{n-1} + \cdots +c_1) = -c_0 \in A^{\f}
\setminus 0$. By part (a), $a \in A^{\f}$.\\
(c): Use the iterative property of $D$ to check that
$D^j(D^i(a))=0$ whenever $j > \deg_{\f}(a) - i$.
\end{proof}

\begin{lemma}\label{L:facts2}
Let $\f$ be a nontrivial exponential map (i.e not the standard
inclusion) on a domain $A$ over $\K$ with $\ch(\K)=p \geq 0$. Let
$x \in A$ have minimal positive $\f$-degree $n$.

\renewcommand{\theenumi}{\alph{enumi}}
\begin{itemize}
\item[(a)] $D^i(x) \in A^{\f}$ for each $i \geq 1$. Moreover,
$D^i(x)=0$ whenever $i > 1$ is not a power of $p$.\label{facts2a}

\item[(b)] If $a \in A \setminus 0$, then $n$ divides
$\deg_{\f}(a)$.\label{facts2b}

\item[(c)] Let $c = D^n(x)$. Then $A$ is a subalgebra of
$A^{\f}[c^{-1}][x]$, where $A^{\f}[c^{-1}] \subseteq
\Frac(A^{\f})$ is the localization of $A^{\f}$ at
$c$.\label{facts2c}

\item[(d)] Let $\trdeg_{\K}$ denote transcendence degree over
$\K$. If $\trdeg_{\K}(A)$ is finite, then $\trdeg_{\K}(A^{\f}) =
\trdeg_{\K}(A) - 1$.\label{facts2d}
\end{itemize}
\end{lemma}

\begin{proof}
In proving parts (a) and (b) we will utilize the following fact.
If $p$ is prime and $i = p^j q$ for some natural numbers $i,j,q$,
then $\binom{i}{p^j} \equiv q \pmod{p}$~\cite[Lemma
5.1]{isaacs}.\\
(a): By part (c) of Lemma~\ref{L:facts}, $D^i(x) \in A^{\f}$ for
all $i \geq 1$. If $p=0$ then $n=1$, for given any element in $A
\setminus A^{\f}$ we can find an element with $\f$-degree 1 by
applying the locally nilpotent derivation $D^1$ sufficiently many
times. In this case, the second statement is immediate. Suppose
now that $p$ is prime and that $i>1$ is not a power of $p$, say $i
= p^j q$, where $j$ is a nonnegative integer and $q \geq 2$ is an
integer not divisible by $p$. Then $D^{i - p^j}(x) \in A^{\f}$ and
\[
0 = D^{p^j} D^{i-p^j}(x) = \binom{i}{p^j}D^i(x) = q D^i(x).
\]
We can divide by $q$ to conclude that $D^i(x)=0$.\\
(b): Again if $p=0$ then $n=1$ and the claim is obvious. Assume
that $p$ is prime. By part (a) we have $n = p^m$ for some integer
$m \geq 0$. If $m=0$, the claim is immediate. Assume that $m>0$.
Let $d = \deg_{\f}(a)$. Suppose that $p$ does not divide $d$. By
part (c) of Lemma~\ref{L:facts}, $\deg_{\f}(D^{d-1}(a)) \leq 1$.
Now, $D^1D^{d-1}(a) = d D^d(a) \ne 0$. So $\deg_{\f}(D^{d-1}(a)) =
1 < n$, contradicting the minimality of $n$. Hence we can write $d
= p^k d_1$ with $k \geq 1$ and $d_1$ not divisible by $p$. Making
a similar computation, $D^{p^k} D^{d-p^k}(a) = d_1 D^d(a) \ne 0$.
This implies that $\deg_{\f}(D^{d-p^k}(a)) = p^k$. Since $n = p^m$
is minimal,
we must have $k \geq m$, and so $n$ divides $d$.\\
(c): Let $a \in A \setminus 0$. By part (b) we can write
$\deg_{\f}(a) = ln$ for some natural number $l$. If $l = 0$ then
$a \in A^{\f}$ and we are done. We use induction on $l>0$.
Elements $c^l a$ and $D^{ln}(a) x^l$ both have $\f$-degree $ln$.
Let us check that $D^{ln}(c^l a) = D^{ln}(D^{ln}(a) x^l)$. First,
$D^{ln}(c^l a) = c^l D^{ln}(a)$ by the Leibniz rule and because
$c^l$ is $\f$-invariant. Secondly, since $D^{ln}(x^l) = D^n(x)^l =
c^l$ and $D^{ln}(a)$ is $\f$-invariant, we see that
$D^{ln}(D^{ln}(a) x^l) = c^l D^{ln}(a)$ as well. (Remark: Though
the equality $D^{ln}(x^l) = D^{n}(x)^l$ does follow from the
Leibniz rule, it may be more immediately observed as follows.
$D^{n}(x)$ is the leading $U$-coefficient of $\f(x)$, and $\f$ is
a homomorphism. Hence the leading $U$-coefficient of $\f(x^l)$ is
also that of $\f(x)^l$.) Therefore, the element $y = c^l a -
D^{ln}(a)x^l$ has $\f$-degree less than $ln$ and hence less than
or equal to $(l-1)n$. By the inductive hypothesis, $y \in
A^{\f}[c^{-1}][x]$. So $a = c^{-l}(y
+ D^{ln}(a) x^l) \in A^{\f}[c^{-1}][x]$.\\
(d): This is immediate from part (c), together with part (b) of
Lemma~\ref{L:facts} which states that $A^{\f}$ is algebraically
closed in $A$.
\end{proof}

\subsection*{Weights}

We often produce a degree function on a domain $A$ by assigning
degree values to some specified generators. On a product of
generators the degree is then defined by the above property (i) of
degree functions. Each element of $A$ can be expressed as a
summation of linearly independent terms, each of which is a
product of generators. The degree of such an expression is then
defined to be the highest degree occurring among the terms. This
is the case with the usual degree functions on polynomials, which
are defined by assigning values to the indeterminates. In this
situation we say that the degree function is obtained by assigning
\emph{weights} to some generators. We will use the idea of weights
repeatedly in proving the results of this paper.

\subsection*{Homogenization of an exponential map}

Let $A$ be a domain over $\K$. Let $\mathbf{Z}$ denote the
integers. Suppose that $A$ has a $\mathbf{Z}$-filtration $\{ A_n
\}$. This means that $A$ is the union of linear subspaces $A_n$
with these properties.
\begin{itemize}
\item[(i)] $A_i \subseteq A_j$ whenever $i \leq j$.

\item[(ii)]  $A_i \cdot A_j \subseteq A_{i+j}$ for all $i,j \in
\mathbf{Z}$.

\item[(iii)] $\bigcap_{n \in \mathbf{Z}} A_n = 0$.
\end{itemize}
Additionally, suppose that
\[
(A_i \setminus A_{i-1}) \cdot (A_j \setminus A_{j-1}) \subseteq
A_{i+j} \setminus A_{i+j-1}
\]
for all $i,j \in \mathbf{Z}$. This will be the case if the
filtration is induced by a degree function. Suppose also that
$\chi$ is a set of generators for $A$ over $\K$ with the following
property: if $a \in A_i \setminus A_{i-1}$ then we can write $a =
\sum_{I} c_{I} \mathbf{x}^{I}$, a summation of monomials $c_{I}
\mathbf{x}^{I}$ built from $\chi$ which are all contained in
$A_i$. This is not an unreasonable property. It merely asserts
some homogeneity on the generating set $\chi$.

Given $a \in A \setminus 0$ there exists $i \in \mathbf{Z}$ for
which $a \in A_i \setminus A_{i-1}$. Write
\[
\tp{a} = a + A_{i-1} \in A_i / A_{i-1},
\]
the \emph{top part of $a$}. We can construct a graded $\K$-algebra
\[
\gr(A) = \bigoplus_{n \in \mathbf{Z}} A_n / A_{n-1}.
\]
Addition on $\gr(A)$ is given by its vector space structure. Given
$\tp{a} = a + A_{i-1}$ and $\tp{b} = b + A_{j-1}$, define $\tp{a}
\, \tp{b} = ab + A_{i+j-1}$. Note that $\tp{a} \, \tp{b} =
\tp{ab}$. Extend this multiplication to all of $\gr(A)$ by the
distributive law. By our assumption on the filtration, $\gr(A)$ is
a domain. Also, $\gr(A)$ is generated by the top parts of the
elements of $\chi$. Therefore, if $\chi$ is a finite set then
$\gr(A)$ is an affine domain.

Let $\grdeg$ be the degree function induced by the grading on
$\gr(A)$. By assigning a weight to an indeterminate $U$ -- call
the weight $\grdeg(U)$, we can extend the grading on $\gr(A)$ to
$\gr(A)[U]$. Given an exponential map $\f: A \to A[U]$ on $A$, the
goal is to obtain an exponential map $\tp{\f}$ on $\gr(A)$. For $a
\in A$, let $\grdeg(a)$ denote $\grdeg(\tp{a})$. Note that
$\grdeg(\tp{a}) = i$ if and only if $a \in A_i \setminus A_{i-1}$.
Consequently, $\grdeg$ can also be viewed as a degree function on
$A$ and on $A[U]$ once the value of $\grdeg(U)$ is determined.

Define
\begin{equation}\label{E:grdegU}
\grdeg(U) = \min \left\{ \frac{\grdeg(x)-\grdeg(D^{i}(x))}{i}
\biggm | x \in \chi, i \in \mathbf{Z}^+ \right\}.\tag{$\star$}
\end{equation}
Let us assume now that $\grdeg(U)$ does exist, i.e. is a rational
number. This will indeed occur whenever $\chi$ is a finite set, as
will be the case with the Danielewski surfaces. If $x \in \chi$
and $n$ is a natural number, then $\grdeg(D^n(x)U^n) \leq
\grdeg(x)$ by our choice of $\grdeg(U)$. From this it follows by
straightforward calculation that $\grdeg(D^n(a) U^n) \leq
\grdeg(a)$ for all $a \in A$ and all natural numbers $n$. (Here we
use the homogeneity assumption imposed on $\chi$.) The reader can
easily work out the details or refer to \cite{thesis}. Note that
this inequality is sharp since
\[
\grdeg(U) = \frac{1}{n}(\grdeg(x) - \grdeg(D^n(x)))
\]
for some $x \in \chi$ and some positive integer $n$ (and also
since $D^0(a) = a$ for all $a \in A$).

For $a \in A$, let
\[
S(a) = \{n \, | \, \grdeg(D^n(a)) + n \grdeg(U) = \grdeg(a) \}.
\]
Define
\[
\tp{\f} (\tp{a}) = \sum_{n \in S(a)} \tp{D^n(a)} U^n
\]
and extend this linearly to define $\tp{\f} : \gr(A) \to
\gr(A)[U]$, the \emph{homogenization\/} or \emph{top part of
$\f$}. One can verify that $\tp{\f}$ is an exponential map on
$\gr(A)$. Refer to~\cite{DHML} for the case $A = \K^{[n]}$. The
proof of the general case is symbolically identical. Let
$\tp{A^{\f}}$ denote the domain generated by the top parts of all
elements in $A^{\f}$. The end result is

\begin{theorem}[H.~Derksen, O.~Hadas, L.~Makar-Limanov~\cite{DHML}]\label{T:toppart}
Let $A$ be a domain over $\K$ with $\mathbf{Z}$-filtration $\{ A_n
\}$ such that $(A_i \setminus A_{i-1}) \cdot (A_j \setminus
A_{j-1}) \subseteq A_{i+j} \setminus A_{i+j-1}$ for all $i,j \in
\mathbf{Z}$. Let $\f$ be a nontrivial exponential map on $A$.
Assume that $\grdeg(U)$ exists as defined above. Then $\tp{\f}$ as
defined above is a nontrivial exponential map on $\gr(A)$.
Moreover, $\tp{A^{\f}}$ is contained in $\gr(A)^{\tp{\f}}$.
\end{theorem}

An important special case of homogenization is when $A$ itself is
graded. Then we can filter $A$ so that $\gr(A)$ is canonically
isomorphic to $A$, and we can choose $\chi$ to be a set of
homogeneous generators of $A$. In this case the top part of $\f$
is a nontrivial exponential map on $A$ (assuming $\grdeg(U)$
exists).

\begin{example}\label{E:toppart} Let $A = \K[X,Y]$, where $\ch(\K) = p$, prime.
Define $\f \in \EXP(A)$ by $\f(X) = X$ and $\f(Y) = Y + U + X
U^p$. We can grade $A$ by assigning weights $\grdeg(X)=\alpha$ and
$\grdeg(Y)=\beta$ (with $\grdeg(\lambda)=0$ for all $\lambda \in
\K^*$, and $\grdeg(0)= -\infty$). Since $\grdeg(D^i(X))= - \infty$
for all $i \geq 1$, $X$ will not contribute to the value of
$\grdeg(U)$. Therefore,

\begin{align*}
\grdeg(U) &= \min \left\{ \frac{\grdeg(Y)-\grdeg(1)}{1},
\frac{\grdeg(Y)-\grdeg(X)}{p} \right\}\\
          &= \min \left\{ \beta, \frac{\beta - \alpha}{p}
          \right\}.
\end{align*}
In any case, $\tp{\f}(X) = X$. If $\beta < \frac{1}{p}(\beta -
\alpha)$ then $\grdeg(U) = \beta$ and $\tp{\f}(Y) = Y + U$. If
$\beta = \frac{1}{p}(\beta - \alpha)$ then $\grdeg(U) = \beta$ and
$\tp{\f}(Y) = \f(Y)$. If $\beta > \frac{1}{p}(\beta - \alpha)$
then $\grdeg(U) = \frac{1}{p}(\beta - \alpha)$ and $\tp{\f}(Y) = Y
+ X U^p$.
\end{example}

\section{Exponential maps of $x^n y = z^2 + h(x)z$}

Let $\K$ be a field with characteristic $p \geq 0$. Let
\[
R = \K[X,Y,Z] / (X^n Y - Z^2 - h(X)Z),
\]
where $n \geq 2$ and $h(X) \in \K[X]$ with $h(0) \ne 0$. Assume
that $\deg_X(h(X)) < n$. (If $d = \deg_X(h(X)) \geq n$ then we can
replace $Y$ by $Y + h_0 X^{d-n}Z$, where $h_0$ is the leading
coefficient of $h(X)$, and replace $h(X)$ by $h(X) - h_0 X^d$.
Iterating this process finitely many times, we can replace $h(X)$
by a polynomial with $X$-degree smaller than $n$.) Let $x,y,z \in
R$ denote the cosets of $X,Y,Z$, respectively. To study the
exponential maps of $R$, we will use filtrations and
homogenization of exponential maps.

\begin{theorem}\label{T:exp} If $\f: R \to R[U]$
is a nontrivial exponential map on $R$, then $R^{\f} =\K[x]$.
Moreover, $z$ has minimal positive $\f$-degree and $D^i(z)$ is
divisible by $x^n$ for each $i \geq 1$, so that
\[
\f(z) = z + x^n f_1(x)U + \sum_{i=1}^{j} x^n f_{p^i}(x) U^{p^i}
\]
for some polynomials $f_i$ and some $j \geq 1$. In fact, any
choice of $j,f_1,f_p,\ldots,f_{p^j}$ in this formula determines a
nontrivial exponential map on $R$.
\end{theorem}

\begin{corollary}
$\ak(R) = \K[x]$.
\end{corollary}

\begin{proof}[Proof of Theorem~\ref{T:exp}]
One easily verifies the final sentence of the theorem by applying
$\f$ to the relation $x^n y = z^2 + h(x)z$ (with $\f(x)=x$),
solving for $\f(y) \in R$, and checking the exponential properties
on the generators $x,y,z$.

Suppose $\f: R \to R[U]$ is a nontrivial exponential map on $R$.
We know that $R^{\f}$ is a subalgebra of $R$ with transcendence
degree 1 over $\K$. So in order to show that $R^{\f} = \K[x]$ it
suffices to show that $x \in R^{\f}$.

View $R$ as a subalgebra of $\K[x,x^{-1},z]$ with $y = x^{-n}(z^2
+ h(x)z)$. Introduce a degree function $w_1$ given by the weights
$w_1(x)=0$, $w_1(y)=2$, and $w_1(z)=1$, with $w_1(\lambda) = 0$
for all $\lambda \in \K^*$ and $w_1(0) = -\infty$, and consider
the $\mathbf{Z}$-filtration $\{R_i\}$ on $R$ induced by these
weights, namely $R_i = \{r \in R \, | \, w_1(r) \leq i \}$.
Observe that $\tp{y} = \tp{x}^{-n}\tp{z}^2$. The graded domain
$\gr(R)$ which corresponds to $w_1$ is generated by $\tp{x},
\tp{y}, \tp{z}$ and subject to the relation $\tp{x}^n \tp{y} =
\tp{z}^2$. Writing $x,y,z$ in place of $\tp{x}, \tp{y}, \tp{z}$,
respectively, we have
\[
\gr(R) = \K[x,y,z] / (x^n y - z^2).
\]

\begin{sublemma}
$R^{\f} \subset \K[x,z]$.
\end{sublemma}

\begin{proof}
Suppose not. Let $f \in R^{\f}$ such that $f \notin \K[x,z]$.
Since $z^2 = x^n y - h(x)z$ we can write $f = f_1(x,y) + z
f_2(x,y)$ for some polynomials $f_1,f_2$. Now $\tp{f}$ is either
$\tp{f_1(x,y)}$ or $\tp{z f_2(x,y)}$, since $f_1(x,y)$ has even
weight while $z f_2(x,y)$ has odd weight. Therefore, $\tp{f} = y^i
z^j g(x)$, where $i$ is a positive integer, $j$ is 0 or 1, and $g$
is some polynomial. (The number $i$ cannot be 0 because $y$
carries the heaviest weight, and our assumption on $f$ is that
some term of either $f_1(x,y)$ or $f_2(x,y)$ must involve $y$.)
Since $R$ is finitely generated by $\chi = \{x,y,z\}$, it is clear
that the value $\grdeg(U)$ exists as defined by formula
(\ref{E:grdegU}). By Theorem~\ref{T:toppart}, the map $\f$ induces
a nontrivial exponential map $\tp{\f}$ on $\gr(R)$ with $\tp{f}
\in \gr(R)^{\tp{\f}}$. Since all factors of $\tp{f}$ belong to
$\gr(R)^{\tp{\f}}$, it follows that $y \in \gr(R)^{\tp{\f}}$.
Since $\gr(R)^{\tp{\f}}$ has transcendence degree 1 over $\K$, we
must have $\gr(R)^{\tp{\f}} = \K[y]$.

Suppose $r$ is an arbitrary element of $R$. Just as with the
element $f$ above we have $\tp{r} = y^i z^j g(x)$ for some natural
numbers $i,j$ and some polynomial $g$. Introduce a new grading on
$\gr(R)$ by the weights $w_2(x)=-1$, $w_2(y) =n$, and $w_2(z) =
0$. This gives us a graded domain $\gr(\gr(R))$ which is naturally
isomorphic to $\gr(R)$. Let us write $\gr(R)$ in place of
$\gr(\gr(R))$, and let us continue to write $x,y,z$ in place of
$\tp{x},\tp{y},\tp{z}$. Under these new weights the top part of
element $\tp{r}$ is $\tp{\tp{r}} = \lambda x^i y^j z^k$ for some
natural numbers $i,j,k$ and some $\lambda \in \K$. The effect of
imposing $w_2$ on $\gr(R)$ is to refine the top parts that were
obtained via $w_1$ in such a way that the top part of every
element from $R$ is a monomial in $x,y,z$. The reader should be
mindful that now on the elements $x,y,z$ there are two weights:
the primary weights given by $w_1$ and the secondary weights given
by $w_2$. By Theorem~\ref{T:toppart}, we obtain a new exponential
map $\tp{\tp{\f}}$ on $\gr(R)$, a ``refinement'' of $\tp{\f}$. But
let us avoid this double-bar notation and use $\tp{\f}$ to denote
the homogenization of $\f$ under the primary and secondary weights
$w_1$ and $w_2$.

Let $a = \deg_{\tp{\f}}(x)$ and $b = \deg_{\tp{\f}}(z)$. Let $D =
\{D^m\}$ be the locally finite iterative higher derivation
associated to $\tp{\f}$. Since $\deg_{\tp{\f}}(y)=0$, the relation
$x^n y = z^2$ indicates that $an = 2b$. Applying $\tp{\f}$ to this
relation and examining the highest power of $U$ which appears
(that being $U^{an} = U^{2b}$), we see that
\[
(D^a(x))^n y = (D^b(z))^2.
\]
Now $D^a(x), D^b(z) \in \gr(R)^{\tp{\f}} = \K[y]$. Also, both
$D^a(x), D^b(z)$ are top parts of elements of $R$ by the way in
which $\tp{\f}$ is defined. Thus each side of the above equation
is a monomial in $y$. If $n$ is even, then the left side has odd
$y$-degree while the right side has even $y$-degree, bringing us
to a contradiction.

Assume that $n$ is odd. We will now argue that $x$ must have
minimal positive $\tp{\f}$-degree, and then we will bring this to
a contradiction. The $\tp{\f}$-degree of $x$ is $a$. To show that
no element has positive $\tp{\f}$-degree smaller than $a$, we will
show that $D^l$ is identically zero for $1 \leq l < a$. Because
$\gr(R)$ is generated by $x,y,z$, it suffices to check $D^l$ on
these elements. Of course $D^l(y) =0$ for all $l \geq 1$. It
remains to study $x$ and $z$.

Write $D^a(x) = \lambda y^i$ for some $\lambda \in \K^*$ and some
natural number $i$. From the above equation we see that $D^b(z) =
\lambda^{n/2} y^{(in + 1)/2}$. Also, since $an = 2b$ we know that
$2$ divides $a$ and $n$ divides $b$. Write $a =2k$. Then $b =nk$.
By the homogeneity of $\tp{\f}$ under both weights $w_1$ and
$w_2$, we know that
\[
w_{\tau}(x) = w_{\tau}(D^a(x)U^a)
\]
for $\tau = 1,2$. This can be rewritten as
\[
w_{\tau}(x) = i w_{\tau}(y) + 2k w_{\tau}(U)
\]
for $\tau =1,2$, where $w_{\tau}(U)$ represents the value
$\grdeg(U)$ for the appropriate grading.

Recall that $w_1(x)=0$ and $w_1(y)=2$. So with $\tau =1$ the above
equations indicate that
\[
w_1(U) = - \frac{i}{k}.
\]
Recall that $w_2(x)=-1$ and $w_2(y)=n$, and so with $\tau = 2$ we
obtain
\[
w_2(U) = - \frac{in + 1}{2k}.
\]

Suppose now that $D^l(x) \ne 0$ for some positive integer $l$. The
homogeneity of $\tp{\f}$ again means that
\[
w_{\tau}(x) = w_{\tau}(D^l(x)U^l)
\]
for $\tau = 1,2$. Also, $D^l(x)$ is a monomial, because by the
definition of $\tp{\f}$ it is the top part of some element from
$R$. Write $D^l(x) = \mu x^\alpha y^\beta z^\gamma$ for some $\mu
\in \K^*$ and some natural numbers $\alpha, \beta, \gamma$. In
fact, $\alpha$ and $\gamma$ must be zero by part (c) of
Lemma~\ref{L:facts} which states that $\deg_{\tp{\f}}(D^l(x)) \leq
2k -l$. We can therefore write
\[
w_{\tau}(x) = \beta w_{\tau}(y) + l w_{\tau}(U)
\]
for $\tau = 1,2$, which in turn becomes the system
\begin{align*}
0 &= 2 \beta - \frac{il}{k}\\
-1 &= n \beta - \frac{iln + l}{2k}.
\end{align*}
We solve this system for $l$ to obtain $l = 2k$. Thus $D^l(x) = 0$
when $1 \leq l < 2k$.

We now argue similarly with $z$ in place of $x$. Suppose that
$D^l(z) \ne 0$ for some positive integer $l$. Again $D^l(z)$ is a
monomial by the homogeneity of $\tp{\f}$. We can write $D^l(z) =
\mu x^{\alpha} y^{\beta}$ for some $\mu \in \K^*$ and some natural
numbers $\alpha, \beta$. ($z$ cannot appear as a factor in
$D^l(z)$ again by part (c) of Lemma~\ref{L:facts}.) Applying $w_1$
and $w_2$ to this presentation of $D^l(z)$ yields the system
\begin{align*}
1 &= 2 \beta - il\\
0 &= - \alpha + n \beta - \frac{iln + l}{2k}.
\end{align*}
We solve this system for $l$ to obtain $l = k(n-2 \alpha)$.
Therefore $l$ must be an odd multiple of $k$. In particular, we
see that $D^l(z) = 0$ for $1 \leq l < k$ and $k< l \leq 2k$.

Let us briefly consider two cases.  First, suppose $\ch(\K) = p
\ne 2$. We already saw that $D^l(x) = 0$ for $1 \leq l < 2k$, and
so applying $D^k$ to the relation $x^n y = z^2$ yields $0 = 2z
D^k(z)$. Hence $D^k(z) = 0$. We now see that $D^l(x)$, $D^l(y)$,
and $D^l(z)$ are identically zero for $1 \leq l < 2k$. Thus $D^l$
is identically zero in that range. As previously discussed, this
means $x$ must be an element of minimal positive $\tp{\f}$-degree
when $p \ne 2$. We now obtain the same conclusion for $p=2$.
Suppose that there does exist an element in $\gr(R)$ of minimal
positive $\tp{\f}$-degree smaller than $2k$. By our analysis of
$D^l(x)$ and $D^l(z)$, that element necessarily must have
$\tp{\f}$-degree $k$. By part (a) of Lemma~\ref{L:facts}, $k$ must
be a power of 2. Also recall that $n$ is odd. Consequently,
$\binom{nk}{k} \equiv n \equiv 1 \pmod{2}$. (We used this fact
about binomial coefficients previously. Refer to the beginning of
the proof of Lemma~\ref{L:facts}.) Now, since $n-1$ is even we
know that $D^{(n-1)k}(z) = 0$. Thus
\[
0 = D^k D^{(n-1)k}(z) = \binom{nk}{k} D^{nk}(z) = D^{nk}(z) \ne 0.
\]
With this contradiction, we now conclude that $x$ must be an
element of minimal positive $\tp{\f}$-degree $2k$ for arbitrary
characteristic $p$.

By part (b) of Lemma~\ref{L:facts2}, we then see that $2k$ must
divide the $\tp{\f}$-degree of $z$, namely $nk$. But this implies
that $2$ divides $n$, contradicting our assumption that $n$ is
odd. This contradiction finishes the proof of the sublemma.
\end{proof}

Let us continue with the proof of Theorem~\ref{T:exp}. Suppose now
that $f \in R^{\f}$ but $f \notin \K[x]$. Since $z^2 = x^n y -
h(x)z$ and $R^{\f} \subset \K[x,z]$, we can write $f = f_1(x) + z
f_2(x)$ for some polynomials $f_1, f_2$ with $f_2 \ne 0$. Once
again we consider $\gr(R)$ induced by $w_1$ and the nontrivial
exponential map $\tp{\f}$ induced by $\f$. Since $w_1(x) = 0$ and
$w_1(z) = 1$, we have $\tp{f} = z f_2(x)$. Since
$\gr(R)^{\tp{\f}}$ is factorially closed and $\tp{f} \in
\gr(R)^{\tp{\f}}$, we must have $z \in \gr(R)^{\tp{\f}}$. Thus
$x^n y = z^2 \in \gr(R)^{\tp{\f}}$, and this implies that $x,y \in
\gr(R)^{\tp{\f}}$. But this means $\tp{\f}$ is trivial, a
contradiction. Therefore $R^{\f}$ is contained in $\K[x]$. Since
$R^{\f}$ is algebraically closed in $R$, we see that $x \in
R^{\f}$ and $R^{\f} = \K[x]$.

Now let us check that $z$ is an element of minimal positive
$\f$-degree. Let $s \in R$ have minimal positive degree. By part
(c) of Lemma~\ref{L:facts2}, there exists $c \in \K[x]$ such that
$R \subseteq \K[x][c^{-1}][s]$. So $R \subseteq \K(x)[s]$ and in
particular $z \in \K(x)[s]$. On the other hand, viewing $y =
x^{-n}(z^2 + h(x)z)$ we know that $R \subseteq \K(x)[z]$, and thus
$s \in \K(x)[z]$. This implies that $z = as + b$ for some $a, b
\in \K(x)$, and since $\deg_{\f}(x)=0$ we have $\deg_{\f}(z) =
\deg_{\f}(s)$.

Since $z$ has minimal positive $\f$-degree, $D^i(z) \in \K[x]$ for
all $i \geq 1$ by part (a) of Lemma~\ref{L:facts2}. If $k \geq 1$,
then
\begin{align*}
x^n D^k(y) &= \sum_{i=0}^k D^i(z)D^{k-i}(z) + h(x)D^k(z)\\
            &= 2z D^k(z) + \left(\sum_{i=1}^{k-1}D^i(z)D^{k-i}(z)
            + h(x)D^k(z)\right).
\end{align*}
So the right hand side of the above equation is a linear (or
possibly constant) polynomial in $z$ with coefficients in $\K[x]$,
and both of these coefficients must be divisible by $x^n$. Recall
that $h(0) \ne 0$ by assumption, so $x$ does not divide $h(x)$. By
induction on $k$ we see that each $D^k(z)$ is divisible by $x^n$.
We have checked for each $i \geq 1$ that $D^i(x) = x^n f_i(x)$ for
some polynomial $f_i$. By part (a) of Lemma~\ref{L:facts2}, $f_i =
0$ whenever $i \geq 2$ is not a power of $p$. Thus we obtain the
formula for $\f(z)$ given in the statement of the theorem.
\end{proof}

Among the exponential maps on $R$ as described by
Theorem~\ref{T:exp} we have those given by $x \mapsto x$ and $z
\mapsto z + x^n f(x)U$. When $\ch(\K)=0$ these are all the
exponential maps on $R$. When $\ch(\K)$ is prime, we can rewrite
the formula for $\f(z)$ in the statement of Theorem~\ref{T:exp} as
\[
z \mapsto z + x^n \gcd(f_i) P(x,U)
\]
for a polynomial $P(x,U)$ which can be viewed as a new
indeterminate $V$. Of course this change of variables will take
several different exponential maps to the same new map.

\section{Automorphisms of $x^n y = z^2 + h(x)z$}

As in the previous section, let $R = \K[X,Y,Z] / (X^n Y - Z^2 -
h(X)Z)$, where $\K$ is a field with characteristic $p \geq 0$, $n
\geq 2$, and $h(X) \in \K[X]$ with $h(0) \ne 0$ and $\deg_X(h(X))
< n$. Let $x,y,z \in R$ denote the cosets of $X,Y,Z$,
respectively. The objective of this section is to describe the
group $\Aut(R)$ of $\K$-algebra automorphisms of $R$ using the
results of the previous section. Let us begin with

\begin{lemma}\label{L:autlemma}
Let $\alpha \in \Aut(R)$. Then $\alpha(x)=\mu x$ for some $\mu \in
\K^*$ such that $h(\mu x)=h(x)$, and either
\begin{itemize}
\item[(a)] $\alpha(z) = z + f(x)$ for some $f \in \K[x]$ with
$f(x) \equiv 0 \pmod{x^n}$, or

\item[(b)] $\alpha(z) = -z + f(x)$ for some $f \in \K[x]$ with
$f(x) \equiv -h(x) \pmod{x^n}$.
\end{itemize}
\end{lemma}

\begin{proof}
If $\f = \sum_i U^i D^i$ is an exponential map on $R$, then
$\alpha^{-1} \f \alpha = \sum_i U^i \alpha^{-1}D^i \alpha$ is
again an exponential map on $R$. Note that if $r \in R$, then the
$(\alpha^{-1} \f \alpha)$-degree of $r$ is equal to the
$\f$-degree of $\alpha(r)$. From this it follows that $\alpha(z)$
must be an element of minimal positive $\f$-degree. Indeed, if $s
\in R$ has lower positive $\f$-degree than that of $\alpha(z)$,
then
\begin{align*}
\deg_{\alpha^{-1} \f \alpha}(\alpha^{-1}(s)) &= \deg_{\f}(s)\\
                            &< \deg_{\f}(\alpha(z))\\
                            &= \deg_{\alpha^{-1} \f \alpha}(z).
\end{align*}
By Theorem~\ref{T:exp}, $z$ has minimal positive $(\alpha^{-1} \f
\alpha)$-degree. This must mean that $\deg_{\alpha^{-1} \f
\alpha}(\alpha^{-1}(s)) \leq 0$, i.e. $\alpha^{-1}(s)$ is
invariant under $\alpha^{-1} \f \alpha$. But then $\alpha^{-1}(s)
\in \K[x]$. We know that $\alpha$ restricts to an automorphism of
$\ak(R) = \K[x]$. Thus $s \in \K[x] = R^{\f}$, contradictory to
our choice of $s$.

Let us fix $\f = \sum_i U^i D^i$ to be the exponential map on $R$
given by $\f(z) = z + x^n U$. Since $z$ and $\alpha(z)$ are both
elements of minimal positive $\f$-degree 1, we have $\alpha(z) =
\lambda z + f$ for some $\lambda, f \in \K[x,x^{-n}]$ (refer to
part (c) of Lemma~\ref{L:facts2}). Actually, we must have
$\lambda, f \in \K[x]$, since the relation $y = x^{-n}(z^2 +
h(x)z)$ in $R$ does not allow for negative powers of $x$ to appear
in the linear polynomial $\lambda z + f$. Moreover, $\lambda \in
\K^*$ because $\alpha$ is invertible.

Since $\alpha$ restricts to an automorphism of $\K[x]$, we have
$\alpha(x) = \mu x + c$ for some $\mu \in \K^*$ and some $c \in
\K$. Now $\deg_{\alpha^{-1} \f \alpha}(z) = \deg_{\f}(\alpha(z)) =
1$, and
\begin{align*}
(\alpha^{-1} D^1 \alpha)(z) &= \alpha^{-1}D^1(\lambda z + f)\\
                            &= \alpha^{-1}(\lambda x^n)\\
                            &= \lambda \mu^{-n} (x - c)^n.
\end{align*}
At the same time, $(\alpha^{-1} D^1 \alpha)(z)$ is divisible by
$x^n$ by Theorem~\ref{T:exp}. Hence $c=0$ and $\alpha(x) = \mu x$.

We are now finished exploiting $\f$. It remains to study how the
relation $x^n y = z^2 + h(x)z$ imposes the remaining information
in the statement of the lemma. Applying $\alpha$ to that relation
we obtain
\begin{align*}
\mu^n x^n \alpha(y) &= \alpha(z)^2 + h(\mu x) \alpha(z)\\
                    &= \lambda^2(z^2 + h(x)z) + g(x,z),
\end{align*}
where
\[
g(x,z) = (2 \lambda f(x) + \lambda h(\mu x) - \lambda^2 h(x))z +
(f(x)^2 + h(\mu x)f(x)).
\]
So
\begin{align*}
\alpha(y) &= \lambda^2 \mu^{-n} x^{-n}(z^2 + h(x)z) +
\mu^{-n}x^{-n} g(x,z)\\
            &= \lambda^2 \mu^{-n} y + \mu^{-n}x^{-n}g(x,z).
\end{align*}
Hence $x^{-n} g(x,z) \in R$. Now $g(x,z)$ is linear as a
polynomial in $z$ with coefficients in $\K[x]$, so $x^{-n}g(x,z)$
cannot have negative $x$-degree. (Again, remember that negative
powers of $x$ can only appear in $R$ when an expression involves
$z^2$ or higher powers of $z$.) Thus $x^n$ must divide $g(x,z)$,
and this means that $x^n$ must divide each coefficient of $g(x,z)$
in $\K[x]$:
\begin{align}
f(x)^2 + h(\mu x)f(x) &\equiv 0 \pmod{x^n},\label{E1}\\
2 \lambda f(x) + \lambda h(\mu x) - \lambda^2 h(x) &\equiv 0
\pmod{x^n}\label{E2}.
\end{align}
To restate equation (\ref{E1}), we know that $x^n$ divides
$f(x)(f(x) + h(\mu x))$. In the following cases we will
demonstrate that $x^n$ must divide either $f(x)$ or $f(x) + h(\mu
x)$. This piece of information allows us to complete the lemma.

\emph{Case (a).} Suppose $f(0)=0$. Then
\[
f(0) + h(\mu \cdot 0) = h(0) \ne 0,
\]
so $x$ does not divide $f(x) + h(\mu x)$, and according to
equation (\ref{E1}) $x^n$ must divide $f(x)$. By equation
(\ref{E2}),
\[
\lambda h(\mu x) - \lambda^2 h(x) \equiv 0 \pmod{x^n}.
\]
Since $\deg_x(h) < n$, we must have
\[
\lambda h(\mu x) - \lambda^2 h(x) = 0.
\]
Setting $x=0$ in this equation, we see that $\lambda = 1$, and
hence $h(\mu x) = h(x)$. We now have $\alpha(z) = z + f(x)$ and
the conditions of the lemma are satisfied.

\emph{Case (b).} Suppose $f(0) \ne 0$. Then $x$ does not divide
$f(x)$, and so $x^n$ must divide $f(x) + h(\mu x)$ by equation
(\ref{E1}). Subsituting
\[
f(x) \equiv -h(\mu x) \pmod{x^n}
\]
in equation (\ref{E2}), we obtain
\[
-\lambda h(\mu x) - \lambda^2 h(x) \equiv 0 \pmod{x^n}.
\]
So
\[
-\lambda h(\mu x) - \lambda^2 h(x) = 0
\]
since $\deg_x(h) < n$. Setting $x=0$ in this equation yields
$\lambda = -1$, and then $h(\mu x) = h(x)$. So $f(x) \equiv -h(x)
\pmod{x^n}$ and $\alpha(z) = -z + f(x)$ as in part (b) of the
lemma.
\end{proof}

We are now in position to prove

\begin{theorem}\label{T:wilkens1}
The group $\Aut(R)$ is generated by
\begin{itemize}
\item[(a)] the automorphisms $E_f$ given by
\begin{align*}
E_f(x) &=x,\\
E_f(y) &= y + 2f(x)z + x^n f(x)^2 + f(x) h(x),\\
E_f(z) &= z + x^n f(x),
\end{align*}
where $f \in \K^{[1]}$,

\item[(b)] the automorphism $T$ given by
\begin{align*}
T(x) &= x,\\
T(y) &= y,\\
T(z) &= -z - h(x),
\end{align*}

\item[(c)] and, if $h(x) = h_1(x^m)$ for some $h_1 \in \K^{[1]}$
and some $m \in \mathbb{N}$, the linear automorphisms $L_{\mu}$
given by
\begin{align*}
L_{\mu}(x) &=\mu x,\\
L_{\mu}(y) &= \mu^{-n}y,\\
L_{\mu}(z) &=z,
\end{align*}
where $\mu \in \K$ such that $\mu^m = 1$.
\end{itemize}
\end{theorem}

\begin{proof}
The map $E_f$ in (a) is obtained by evaluating $U=1$ in the
exponential map given by $\f(z)= z + x^n f(x) U$. It is easy to
check that all of the maps in (b) and (c) are indeed
automorphisms. The conditions $h(x)=h_1(x^m)$ and $\mu^m =1$ in
(c) describe the only possible way that $h(\mu x) = h(x)$ for $\mu
\ne 1$, as in Lemma~\ref{L:autlemma}. If $m = 1$ then (c)
describes only the identity automorphism, and if $m = 0$ then $h$
is constant and $\mu$ can be any element of $\K^*$. Automorphisms
of the forms $L_{\mu} E_f$ and $L_{\mu} T E_f$ cover all possible
maps described in Lemma~\ref{L:autlemma}.
\end{proof}

The set $L$ of automorphisms $L_{\mu}$ in (c) is an abelian
subgroup of $\Aut(R)$. Let $H$ be the subgroup generated by $L$
and $T$. $H$ is an abelian group, the internal direct product
$\langle T \rangle \times L$. Let $N$ be the set of automorphisms
$E_f$ in (a). $N$ is a normal subgroup of $\Aut(R)$, and $\Aut(R)$
is the semi-direct product $N \rtimes H$.

$N$ is isomorphic to $\K^{[1]}$ as an additive group, and $\langle
T \rangle$ is cyclic of order 2. Turning to the conditions of (c),
we see that $L$ is trivial if $m=1$ and isomorphic to $\K^*$ if
$m=0$. Otherwise $L$ is a cyclic group of order dividing $m$. (If
$\K$ is algebraically closed then $L$ has order $m$.) Let $C_k$
denote the cyclic group of order $k$. To summarize:

\begin{corollary} $\Aut(R) \cong \K^{[1]} \rtimes H$, where $H$ is
as follows.
\begin{itemize}
\item[(a)] If $h$ is a constant polynomial, then $H \cong C_2
\times \K^*$.

\item[(b)] If $h(x) = h_1(x^m)$ for some $m>1$, then $H \cong C_2
\times C_k$ for some factor $k$ of $m$.

\item[(c)] Otherwise (for a ``typical'' polynomial $h$), $H \cong
C_2$.
\end{itemize}
\end{corollary}

\section{Remarks on the cancellation problem}

For $i=1,2$ let $R_i = \K[X,Y,Z]/(X^{n_i} Y - h_i(X)Z)$, where
$n_i \geq 2$ and $h_i(0) \ne 0$. As mentioned in the introduction,
these algebras (when $\K = \mathbf{C}$) are known to be a class of
counterexamples to the cancellation problem. That is, $R_1^{[1]}
\cong R_2^{[1]}$, while in general $R_1 \ncong R_2$. This has been
explained geometrically \cite{danielewski,wilkens}, but let us
briefly provide an algebraic explanation. To show the isomorphism
of polynomial rings, we can try the following approach. Embed
$R_1$ in $R_2[T]$, and then find an exponential map $\f: R_2[T]
\to R_2[T][U]$ with ring of invariants $R_1$ and such that $\f(s)
= s + U$ for some $s$. This will imply that $R_2[T] = R_1[s]$ by
part (c) of Lemma~\ref{L:facts2} (because, in the notation of that
lemma, we will have $c=1$). The element $s$ is commonly called a
\emph{slice}.

Here are the formulae for a special case. Assume that $n_1 < n_2
\leq 2n_1$ and that $h_1(X) = h_2(X) = 1$. $\K$ can be any field.
Let $x_i,y_i,z_i$ denote the cosets of $X,Y,Z$ in $R_i$,
respectively. So for $i=1,2$ we have algebras $R_i$ generated by
$x_i,y_i,z_i$ subject to the relations
\begin{align}
x_{1}^{n_1} y_1 &= z_{1}^2 + z_1,\label{R1}\\
x_{2}^{n_2} y_2 &= z_{2}^2 + z_2.\label{R2}
\end{align}
Embed $R_1$ in $R_2$ by sending $x_1$ to $x_2$, $z_1$ to $z_2$,
and $y_1$ to $x_2^{n_2 - n_1}y_2$. Let $\widetilde{R_1}$ denote
the image of $R_1$ in $R_2$. By Theorem~\ref{T:exp} we have an
exponential map on $\widetilde{R_1}$ defined by sending $x_2$ to
$x_2$ and $z_2$ to $z_2 + x_2^{n_1} T$, where $T$ is the
indeterminate which parameterizes the exponential map. Of course
exponential maps are injective, and the composition of the
embedding $R_1 \hookrightarrow \widetilde{R_1}$ with the
exponential map on $\widetilde{R_1}$ gives us an embedding of
$R_1$ in $R_2[T]$ given by
\begin{align*}
x_1 & \mapsto x_2,\\
z_1 & \mapsto z_2 + x_2^{n_1} T,\\
y_1 & \mapsto x_2^{n_2 - n_1} y_2 + (2 z_2 +1)T + x_2^{n_1} T^2.
\end{align*}
Let us identify $R_1$ with its image under this embedding,
yielding $x_1 = x_2 = x$ and
\begin{align}
z_1 &= z_2 + x^n T\label{R3},\\
y_1 &= x^{n_2 - n_1}y_2 + (2z_1 + 1)T - x^{n_1}T^2.\label{R4}
\end{align}
Relations (\ref{R1}), (\ref{R2}), (\ref{R3}), and (\ref{R4})
completely describe $R_2[T]$. In fact, relations (\ref{R2}) and
(\ref{R3}) are unnecessary since from (\ref{R1}) and (\ref{R4}) we
recover the relation
\[
x^{n_2} y_2 = (z_1 - x^{n_1} T)^2 + (z_1 - x^{n_1} T).
\]
This means $R_2[T]$ is generated by $x,y_1,y_2,z_1,T$ and subject
to the relations (\ref{R1}) and (\ref{R4}). Define $\f: R_2[T] \to
R_2[T][U]$ by $\f(x) = x$, $\f(y_1) = y_1$, and $\f(z_1) = z_1$,
with
\begin{align*}
\f(y_2) &= y_2 + (2z_1 - 2x^{n_1}T + 1)U + x^{n_2}U^2,\\
\f(T) &= T - x^{n_2 - n_1}U.
\end{align*}
One can easily observe that $\f$ is an exponential map on $R_2[T]$
with ring of invariants $R_1$ by checking the exponential
properties on the generators. Moreover one can verify that $\f(s)
= s + U$, where
\[
s = -4 x^{3n_1 - n_2} T^3 + 3x^{2n_1 - n_2}(2z_1 + 1)T^2 +
4x^{n_1}y_2 T + y_2(2z_1 + 1).
\]
Consequently, $R_2[T] = R_1[s]$ by part (c) of
Lemma~\ref{L:facts2}. Because $T$ is also an element of minimal
positive $\f$-degree 1, that same lemma leads us to expect that
$s$ should be a linear polynomial in $T$ with coefficients in
$R_1[x^{n_1 - n_2}]$, and indeed one can check that
\[
s = -x^{n_1 - n_2}(T - y_1(2z_1 + 1)).
\]

To conclude, let us prove the following theorem which shows that
two Danielewski surfaces are in general not isomorphic.
\begin{theorem}
Let $\K$ be a field. For $i=1,2$ let $R_i = \K[X,Y,Z]/(X^{n_i} Y -
h_i(X)Z)$, where $n_i \geq 2$ and $h_i(0) \ne 0$. Let
$x_i,y_i,z_i$ denote the cosets of $X,Y,Z$ in $R_i$, respectively.
Then $R_1 \cong R_2$ if and only if $n_1 = n_2$ and $h_2(x) = \eta
h_1(\mu x)$ for some $\eta, \mu \in \K^*$.
\end{theorem}

\begin{proof}
($\Leftarrow$) The map given by $x_1 \mapsto \mu x_2$, $y_1
\mapsto \eta^{-2} \mu^{-n}y_2$, $z_1 \mapsto \eta^{-1}z_2$ clearly
defines an isomorphism of $R_1$ onto $R_2$, where $n=n_1 =n_2$.

($\Rightarrow$) We can again assume that $\deg_{x_i}(h_i(x_i)) <
n_i$, $i=1,2$. We proceed as in Lemma~\ref{L:autlemma}, and the
reader should refer to the arguments made there. Suppose $\alpha:
R_1 \to R_2$ is an isomorphism. If $\f \in \EXP(R_2)$, then
$\alpha^{-1} \f \alpha \in \EXP(R_1)$. Consequently, as in the
proof of Lemma~\ref{L:autlemma} we can conclude that $\alpha(x_1)
= \mu x_2$ and $\alpha(z_1) = \lambda z_2 + f(x_2)$ for some
$\lambda, \mu \in \K^*$ and some polynomial $f_2$. Applying
$\alpha$ to the relation $x_1^{n_1} y_1 = z_1^2 + h_1(x_1)z_1$ we
see that $x_2^{n_1}\alpha(y_1)$ is a polynomial in $x_2$ and $z_2$
with $z_2$-degree 2. The relation $y_2 = x_2^{-n_2}(z_2^2 +
h_2(x_2)z_2)$ on $R_2$ will not allow us to divide by $x_2^{n_1}$
unless $n_1 \leq n_2$. Repeating this analysis with $\alpha^{-1}$
we also must have $n_2 \leq n_1$, so that $n_1 = n_2 = n$. Next,
just as in the proof of Lemma~\ref{L:autlemma} we obtain
\begin{align*}
f(x_2)^2 + h_1(\mu x_2)f(x_2) & \equiv 0 \pmod{x_2^n},\\
2 \lambda f(x_2) + \lambda h_1(\mu x_2) - \lambda^2 h_2(x_2) &
\equiv 0 \pmod{x_2^n}.
\end{align*}
Recall now that $\deg_{x_i}(h_i(x_i)) < n$ for $i=1,2$. Continuing
as with Lemma~\ref{L:autlemma} we consider two possibilities. If
$f(0) = 0$ then we find that $h_2(x_2) = \lambda^{-1} h_1(\mu
x_2)$, and the conditions of the theorem are satisfied with $\eta
= \lambda^{-1}$. If $f(0) \ne 0$ we conclude that $h_2(x_2) = -
\lambda^{-1} h_1(\mu x_2)$, and again we are done with $\eta = -
\lambda^{-1}$.
\end{proof}


\end{document}